\RequirePackage{etoolbox}
\csdef{input@path}{{style/}{graphics/}}
\documentclass[seceqn,MSNbibl,number,citesort,dvips]{arxbj}
\usepackage{graphicx}

\aid{0}
\volume{21}
\issue{1}
\pubyear{2015}
\firstpage{62}
\lastpage{82}
\doi{10.3150/13-BEJ561} 

\makeatletter
\newtheorem{Lem}{Lemma}[section]
\newtheorem{Theor}{Theorem}[section]
\newremark{Setup}{Setup}
\newcommand{\ny}{n\to\infty}
\newcommand{\R}{\mathbb{R}}
\newcommand{\Pub}{P_{[\ub]}}
\newcommand{\mub}{\boldsymbol{\mu}}
\newcommand{\thetab}{\boldsymbol{\theta}}
\newcommand{\Sigmab}{\boldsymbol{\Sigma}}
\newcommand{\Xb}{\mathbf{X}}
\newcommand{\yb}{\mathbf{y}}
\newcommand{\ub}{\mathbf{u}}
\newcommand{\Ab}{\mathbf{A}}
\newcommand{\bb}{\mathbf{b}}

\newcommand{\Unif}{\operatorname{Unif}}
\newcommand{\Supp}{\operatorname{Supp}}
\makeatother

\begin{document}
\begin{frontmatter}

\title{Nonparametrically consistent depth-based classifiers}

\runtitle{Nonparametrically consistent depth-based classifiers}

\begin{aug}
\author[1]{\inits{D.}\fnms{Davy} \snm{Paindaveine}\corref{}\thanksref{E1}\ead[label=E1,mark]{dpaindav@ulb.ac.be}}
\and
\author[1]{\inits{G.}\fnms{Germain} \snm{Van Bever}\thanksref{E2}\ead[label=E2,mark]{gvbever@ulb.ac.be}}
\address[1]{ECARES \& D\'epartement de Math\'ematique, Universit\'e
libre de Bruxelles, Belgium.\\ \printead{E1,E2}}
\end{aug}

\received{\smonth{6} \syear{2012}}
\revised{\smonth{6} \syear{2013}}

%
\begin{abstract}
We introduce a class of depth-based classification procedures that are
of a nearest-neighbor nature. Depth, after symmetrization, indeed
provides the center-outward ordering that is necessary and sufficient
to define nearest neighbors. Like all their depth-based competitors,
the resulting classifiers are affine-invariant, hence in particular are
insensitive to unit changes. Unlike the former, however, the latter
achieve Bayes consistency under virtually any absolutely continuous
distributions -- a concept we call \emph{nonparametric consistency}, to
stress the difference with the stronger \emph{universal consistency}
of the standard $k$NN classifiers. We investigate the finite-sample
performances of the proposed classifiers through simulations and show
that they outperform affine-invariant nearest-neighbor classifiers
obtained through an obvious standardization construction. We illustrate
the practical value of our classifiers on two real data examples.
Finally, we shortly discuss the possible uses of our depth-based
neighbors in other inference problems.
\end{abstract}

%
\begin{keyword}
\kwd{affine-invariance}
\kwd{classification procedures}
\kwd{nearest neighbors}
\kwd{statistical depth functions}
\kwd{symmetrization}
\end{keyword}

\end{frontmatter}
%
\section{Introduction}
\label{introsec}

The main focus of this work is on the standard classification setup in
which the observation, of the form ($\Xb,Y$), is a random vector
taking values in $\R^d\times\{0,1\}$. 
A classifier is a function $m\dvtx \R^d \to\{0,1\}$ that associates with
any value ${\mathbf{x}}$ a predictor for the corresponding
``class'' $Y$. Denoting by $\mathbb{I}[A]$ the indicator function of
the set $A$, the so-called Bayes classifier, defined through
%
\begin{equation}
\label{Bayes} m_\mathrm{Bayes}({\mathbf{x}})= \mathbb{I} \bigl[
\eta({\mathbf{x}})> 1/2 
\bigr] , \qquad \mbox{with } \eta({\mathbf{x}})=P[Y=1 \mid \Xb= {\mathbf{x}}], 
\end{equation}
is optimal in the sense that it minimizes the probability of
misclassification $P[m(\Xb)\neq Y]$. Under absolute continuity
assumptions, the Bayes rule rewrites
%
\begin{equation}
\label{Bayesdens} m_\mathrm{Bayes}({\mathbf{x}}) =\mathbb{I} \biggl[
\frac{f_1({\mathbf{x}})}{f_0(
{\mathbf{x}})}>\frac{\pi_0}{\pi_1} \biggr],
\end{equation}
where $\pi_j=P[Y=j]$ and $f_j$ denotes the pdf of $\Xb$ conditional
on $[Y=j]$. Of course, empirical classifiers $\hat{m}^{(n)}$ are obtained
from \mbox{i.i.d.} copies $(\Xb_i,Y_i)$, $i=1,\ldots,n$, of $(\Xb
,Y)$, and it is desirable that such classifiers are consistent, in the
sense that, as $n\to\infty$, the probability of misclassification of
$\hat{m}^{(n)}$, conditional on $(\Xb_i,Y_i)$, $i=1,\ldots,n$,
converges in probability to the probability of misclassification of
the Bayes rule.
If this convergence holds irrespective of the distribution of $(\Xb
,Y)$, the consistency is said to be \emph{universal}.

Classically, parametric approaches assume that the conditional
distribution of $\Xb$ given $[Y=j]$ is multinormal with mean $\mub_j$
and covariance matrix $\Sigmab_j$ ($j=0,1$). This gives rise to the
so-called \emph{quadratic discriminant analysis (QDA)} -- or to \emph
{linear discriminant analysis (LDA)} if it is further assumed that
$\Sigmab_0=\Sigmab_1$. It is standard to estimate the parameters
$\mub_j$ and $\Sigmab_j$ $(j=0,1)$ by the corresponding sample means
and empirical covariance matrices, but the use of more robust
estimators was recommended in many works; see, for example, Randles \textit{et~al.} \cite{Randles:1978vq}, He and Fung \cite
{He:2000vs}, Dehon and Croux \cite{Dehon:2001um}, or
Hartikainen and Oja \cite{Hartikainen2006}. Irrespective of
the estimators used, however, these
classifiers fail to be consistent away from the elliptical case.

Denoting by $d_{\Sigmab}({\mathbf{x}},\mub)=((
{\mathbf{x}}-\mub)'\Sigmab^{-1}({\mathbf{x}}-\mub))^{1/2}$ the
Mahalanobis distance between ${\mathbf{x}}$ and $\mub$ in the
metric associated with the symmetric and positive definite matrix
$\Sigmab$, it is well known that the QDA classifier rewrites
%
\begin{equation}
\label{qda} m_\mathrm{QDA}({\mathbf{x}}) = \mathbb{I} \bigl[
d_{\Sigmab_1}({\mathbf{x}},\mub_1) < d_{\Sigmab
_0}({
\mathbf{x}},\mub_0)+C \bigr],
\end{equation}
where the constant $C$ depends on $\Sigmab_0$, $\Sigmab_1$, and $\pi
_0$, hence classifies ${\mathbf{x}}$ into Population $1$ if it
is sufficiently more central in Population 1 than in Population 0
(centrality, in elliptical setups, being therefore measured with
respect to the geometry of the underlying equidensity contours). This
suggests that \emph{statistical depth functions}, that are mappings of
the form ${\mathbf{x}}\mapsto D({\mathbf{x}},P)$
indicating how central ${\mathbf{x}}$ is with respect to a
probability measure $P$ (see Section~\ref{depthsec} for a more precise
definition), are appropriate tools to perform nonparametric
classification. Indeed, denoting by $P_j$ the probability measure
associated with Population $j$ ($j=0,1$), (\ref{qda}) makes it natural
to consider classifiers of the form
\[
m_D({\mathbf{x}}) = \mathbb{I} \bigl[ D({\mathbf{x}},P_1)>D({
\mathbf{x}},P_0) \bigr] ,
\]
based on some fixed statistical depth function $D$. This \emph
{max-depth approach} was first proposed in Liu, Parelius and Singh
\cite{Liu:1999wm} and was
then investigated in Ghosh and Chaudhuri \cite{Ghosh:2005wx}.
Dutta and Ghosh \cite{Dutta:2012ig,Dutta:wt}
considered max-depth classifiers based on
the projection depth and on (an affine-invariant version of) the $L^p$
depth, respectively. Hubert and Van~der Veeken \cite{Hubert:2010gk}
modified the max-depth
approach based on projection depth to better cope with possibly skewed data.

Recently, Li, Cuesta-Albertos and Liu \cite{Lietal2012} proposed the
``Depth vs Depth'' (DD)
classifiers that extend the max-depth ones by constructing appropriate
polynomial separating curves in the DD-plot, that is, in the scatter
plot of the points $(D_0^{(n)}(\Xb_i),D_1^{(n)}(\Xb_i))$, $i=1,\ldots,n$,
where $D_j^{(n)}(\Xb_i)$ refers to the depth of $\Xb_i$ with respect to
the data points coming from Population $j$. Those separating curves are
chosen to minimize the empirical misclassification rate on the training
sample and their polynomial degree $m$ is chosen through
cross-validation. Lange, Mosler and Mozharovskyi \cite{Lange12} defined
modified DD-classifiers that
are computationally efficient and apply in higher dimensions (up to
$d=20$). Other depth-based classifiers were proposed in J{\"o}rnsten \cite{Jornsten:2004fv}, Ghosh and Chaudhuri
\cite{Ghosh:2005ts} and Cui, Lin and Yang \cite{Cui:2008jg}.

Being based on depth, these classifiers are clearly of a nonparametric
nature. An important requirement in nonparametric classification,
however, is that consistency holds as broadly as possible and, in
particular, does not require ``structural'' distributional assumptions.
In that respect, the depth-based classifiers available in the
literature are not so satisfactory, since they are at best consistent
under elliptical distributions only.\footnote{The classifiers from
Dutta and Ghosh \cite{Dutta:wt} are an exception that slightly
extends consistency to
(a subset of) the class of $L_p$-elliptical distributions.}
This restricted-to-ellipticity consistency implies that, as far as
consistency is concerned, the Mahalanobis depth is perfectly sufficient
and is by no means inferior to the ``more nonparametric'' (Tukey
\cite{Tukey:1975ts}) halfspace depth or (Liu \cite{Liu:1990p224}) simplicial
depth, despite the fact that it uninspiringly leads to LDA through the
max-depth approach. Also, even this restricted consistency often
requires estimating densities; see, for example, Dutta and Ghosh \cite{Dutta:2012ig,Dutta:wt}. This is somewhat
undesirable since density and
depth are quite antinomic in spirit (a deepest point may very well be a
point where the density vanishes). Actually, if densities are to be
estimated in the procedure anyway, then it would be more natural to go
for density estimation all the way, that is, to plug density estimators
in (\ref{Bayesdens}).

The poor consistency of the available depth-based classifiers actually
follows from their \emph{global} nature. Zakai and Ritov
\cite{Zakai:2009uq} indeed
proved that any universally consistent classifier needs to be of a
\emph{local} nature. In this paper, we therefore introduce local
depth-based classifiers, that rely on nearest-neighbor ideas (kernel
density techniques should be avoided, since, as mentioned above, depth
and densities are somewhat incompatible). From their nearest-neighbor
nature, they will inherit consistency under very mild conditions, while
from their depth nature, they will inherit affine-invariance and
robustness, two important features in multivariate statistics and in
classification in particular. Identifying nearest neighbors through
depth will be achieved via an original symmetrization construction. The
corresponding depth-based neighborhoods are of a nonparametric nature
and the good finite-sample behavior of the resulting classifiers most
likely results from their data-driven adaptive nature.


The outline of the paper is as follows. In Section~\ref
{depthneighsec}, we first recall the concept of statistical depth
functions (Section~\ref{depthsec}) and then describe our
symmetrization construction that allows to define the depth-based
neighbors to be used later for classification purposes (Section~\ref
{neighsec}). In Section~\ref{classisec}, we define the proposed
depth-based nearest-neighbor classifiers and present some of their
basic properties (Section~\ref{definprocsec}) before providing
consistency results (Section~\ref{consistencysec}). In Section~\ref
{simusec}, Monte Carlo simulations are used to compare the
finite-sample performances of our classifiers with those of their
competitors. In Section~\ref{realexsec}, we show the practical value
of the proposed classifiers on two real-data examples. We then discuss
in Section~\ref{finalsec} some further applications of our depth-based
neighborhoods. Finally, the \hyperref[appsec]{Appendix} collects the technical proofs.

\section{Depth-based neighbors}
\label{depthneighsec}

In this section, we review the concept of statistical depth functions
and define the depth-based neighborhoods on which the proposed
nearest-neighbor classifiers will be based.

\subsection{Statistical depth functions}
\label{depthsec}

Statistical depth functions allow to measure \emph{centrality} of any
${\mathbf{x}}\in\R^d$ with respect to a probability measure
$P$ over $\R^d$ (the larger the depth of ${\mathbf{x}}$, the
more central ${\mathbf{x}}$ is with respect to $P$). Following
Zuo and Serfling \cite{Zuo:2000p160}, we define a
statistical depth function as a
bounded mapping $D( \cdot,P)$ from $\R^d$ to $\R^+$
that satisfies the following four properties:
\begin{enumerate}[(P3)]
\item[(P1)] \emph{affine-invariance}: for any $d\times d$ invertible
matrix $\Ab$, any $d$-vector $\bb$ and any distribution $P$ over $\R
^d$, $D(\Ab{\mathbf{x}}+\bb,P^{\Ab,\bb})=D(
{\mathbf{x}},P)$, where $P^{\Ab,\bb}$ is defined through $P^{\Ab
,\bb
}[B]=P[\Ab^{-1}(B-\bb)]$ for any $d$-dimensional Borel set $B$;
\label{Ps}
\item[(P2)] \emph{maximality at center}: for any $P$ that is
symmetric about $\thetab$ (in the sense\footnote{Zuo and Serfling \cite{Zuo:2000p160}
also considers more general symmetry concepts; however, we restrict in
the sequel to central symmetry, that will be the right concept for our
purposes.} that $P[\thetab+B]=P[\thetab-B]$ for any $d$-dimensional
Borel set $B$), $D(\thetab,P)=\sup_{{\mathbf{x}}\in\R^d}
D({\mathbf{x}},P)$;
%
\item[(P3)] \emph{monotonicity relative to the deepest point}: for
any $P$ having deepest point $\thetab$, for any ${\mathbf{x}}\in\R
^d$ and any $\lambda\in[0,1]$, $D({\mathbf{x}},P)\leq D((1-\lambda)
\thetab+\lambda{\mathbf{x}},P)$;
\item[(P4)] \emph{vanishing at infinity}: for any $P$, $D(
{\mathbf{x}},P)\to0$ as $\|{\mathbf{x}}\|\to\infty$.
\end{enumerate}

For any statistical depth function and any $\alpha>0$, the set
$R_\alpha(P)=\{ {\mathbf{x}}\in\R^d: D({\mathbf{x}},P)\geq\alpha
\}$ is called \emph{the depth region of order}
$\alpha$. These regions are nested, and, clearly, inner regions
collect points with larger depth. Below, it will often be convenient to
rather index these regions by their probability content: for any
$\beta\in[0,1)$, we will denote by
$
R^{\beta}(P)
$
the smallest $R_{\alpha}(P)$ that has $P$-probability larger than or
equal to $\beta$.
Throughout, subscripts and superscripts for depth regions are used for
depth levels and probability contents, respectively. 

Celebrated instances of statistical depth functions include
\begin{enumerate}[(iii)]
\item[(i)] the Tukey \cite{Tukey:1975ts}
halfspace depth
$
D_H({\mathbf{x}},P)
=
\inf_{\ub\in\mathcal{S}^{d-1}} P[\ub'(\Xb-{\mathbf{x}})\geq0 ]$,
where $\mathcal{S}^{d-1}=\{\ub\in\R^d:\|\ub\|=1\}$ is the unit
sphere in $\R^d$;
%
\item[(ii)] the Liu \cite{Liu:1990p224}
simplicial depth
$
D_S({\mathbf{x}},P)
=
P[{\mathbf{x}} \in S( \Xb_{1}, \Xb_{2},\ldots, \Xb_{d+1}) ]$,
where $S( {\mathbf{x}}_{1},\allowbreak  {\mathbf{x}}_{2},\ldots,
{\mathbf{x}}_{d+1})$ denotes the closed simplex with vertices
${\mathbf{x}}_{1}, {\mathbf{x}}_{2},\ldots,
{\mathbf{x}}_{{d+1}}$ and where $\Xb_{1}, \Xb_{2},\ldots,
\Xb_{d+1}$ are \mbox{i.i.d.} $P$;
\item[(iii)] the Mahalanobis depth
$
D_M({\mathbf{x}},P)
=
1/(1+d^2_{\Sigmab(P)}({\mathbf{x}},\mub(P)))$,
for some affine-equivariant location and scatter functionals $\mub(P)$
and $\Sigmab(P)$;
\item[(iv)] the projection depth
$
D_{Pr}({\mathbf{x}},P)
=
1/(1+\sup_{\ub\in\mathcal{S}^{d-1}} |\ub'{\mathbf{x}}-\mu
(\Pub)|/\sigma(\Pub))$,
where $\Pub$ denotes the probability distribution of $\ub'\Xb$ when
$\Xb\sim P$ and where $\mu(P)$ and $\sigma(P)$ are univariate
location and scale functionals, respectively.
\end{enumerate}
Other depth functions are the simplicial volume depth, the spatial
depth, the $L_p$ depth, etc. Of course, not all such depths fulfill
properties (P1)--(P4) for any distribution $P$; see Zuo and Serfling \cite{Zuo:2000p160}. A further concept of depth, of a
slightly different
($L_2$) nature, is the so-called \emph{zonoid depth}; see Koshevoy and Mosler \cite{Koshevoy:1997tm}.


Of course, if $d$-variate observations $\Xb_1,\ldots,\Xb_n$ are
available, then sample versions of the depths above are simply obtained
by replacing $P$ with the corresponding empirical distribution
$P^{(n)}$ (the sample simplicial depth then has a $U$-statistic structure).

%

A crucial fact for our purposes is that a sample depth provides a \emph
{center-outward ordering} of the observations with respect to the
corresponding deepest point $\hat{\thetab}^{(n)}$: one may indeed order
the $\Xb_i$'s in such a way that
%
\begin{equation}
\label{s} D\bigl(\Xb_{(1)},P^{(n)}\bigr) \geq D\bigl(
\Xb_{(2)},P^{(n)}\bigr) \geq \cdots \geq D\bigl(
\Xb_{(n)},P^{(n)}\bigr).
\end{equation}
Neglecting possible ties, this states that, in the depth sense, $\Xb
_{(1)}$ is the observation closest to $\hat{\thetab}^{(n)}$, $\Xb_{(2)}$
the second closest, \ldots, and $\Xb_{(n)}$ the one farthest away
from $\hat{\thetab}^{(n)}$.

For most classical depths, there may be infinitely many deepest points,
that form a convex region in $\R^d$. This will not be an issue in this
work, since the symmetrization construction we will introduce, jointly
with properties (Q2)--(Q3) below, asymptotically guarantees unicity of
the deepest point. For some particular depth functions, unicity may
even hold for finite samples: for instance, in the case of halfspace
depth, it follows from Rousseeuw and Struyf \cite
{Rousseeuw:2004p253} and results on the
uniqueness of the symmetry center (Serfling \cite
{Serfling:2006up}) that, under
the assumption that the parent distribution admits a density,
symmetrization implies almost sure unicity of the deepest point.

\subsection{Depth-based neighborhoods}
\label{neighsec}

A statistical depth function, through (\ref{s}), can be used to define
neighbors of the deepest point $\hat{\thetab}^{(n)}$. Implementing a
nearest-neighbor classifier, however, requires defining neighbors of
any point ${\mathbf{x}}\in\R^d$. Property (P2) provides the
key to the construction of an ${\mathbf{x}}$-outward ordering
of the observations, hence to the definition of depth-based neighbors
of ${\mathbf{x}}$: symmetrization with respect to
${\mathbf{x}}$.

More precisely, we propose to consider depth with respect to the
empirical distribution $P^{(n)}_{{\mathbf{x}}}$ associated
with the sample obtained by adding to the original observations $\Xb
_1,\Xb_2,\ldots,\Xb_n$ their reflections $2{\mathbf{x}}-\Xb
_1,\ldots,2{\mathbf{x}}-\Xb_n$ with respect to $
{\mathbf{x}}$. Property (P2) implies that ${\mathbf{x}}$ is
the -- unique (at least asymptotically; see above) -- deepest point with
respect to $P^{(n)}_{{\mathbf{x}}}$. Consequently, this
symmetrization construction, parallel to (\ref{s}), leads to an
(${\mathbf{x}}$-outward) ordering of the form
\[
D\bigl(\Xb_{{\mathbf{x}},(1)},P^{(n)}_{{\mathbf{x}}}\bigr) \geq D\bigl(
\Xb_{{\mathbf{x}},(2)},P^{(n)}_{{\mathbf{x}}}\bigr) \geq \cdots \geq D\bigl(
\Xb_{{\mathbf{x}},(n)},P^{(n)}_{{\mathbf{x}}}\bigr).
\]
Note that the reflected observations are only used to define the
ordering but are not ordered themselves. For any $k\in\{1,\ldots,n\}
$, this allows to identify -- up to possible ties -- the $k$ nearest
neighbors $\Xb_{{\mathbf{x}},(i)}$, $i=1,\ldots,k$, of
${\mathbf{x}}$. In the univariate case ($d=1$), these $k$
neighbors coincide -- irrespective of the statistical depth function
$D$ -- with the $k$ data points minimizing the usual distances
$|X_{i}-x|$, $i=1,\ldots,n$.

In the sequel, the corresponding \emph{depth-based
neighborhoods} -- that is, the sample depth regions $R^{(n)}_{
{\mathbf{x}},\alpha}=R_{\alpha}(P^{(n)}_{{\mathbf{x}}})$ -- will
play an important role. In accordance with the notation from the
previous section, we will write
$
R_{{\mathbf{x}}}^{\beta(n)}
$
for the smallest depth region $R^{(n)}_{{\mathbf{x}},\alpha}$
that contains at least a proportion $\beta$ of the data points $\Xb
_1,\Xb_2,\ldots,\Xb_n$. For $\beta=k/n$, $R_{{\mathbf{x}}}^{\beta
(n)}$ is therefore the smallest depth-based neighborhood
that contains $k$ of the $\Xb_i$'s; ties may imply that the number of
data points in this neighborhood, $K_{\mathbf{x}}^{\beta(n)}$
say, is strictly larger than $k$.

Note that a distance (or pseudo-distance) $({\mathbf{x}},\yb
)\mapsto d({\mathbf{x}},\yb)$ that is symmetric in its
arguments is not needed to identify nearest neighbors of $
{\mathbf{x}}$. For that purpose, a collection of ``distances'' $\yb
\mapsto
d_{\mathbf{x}}(\yb)$ from a fixed point is indeed sufficient
(in particular, it is irrelevant that this distance satisfies or not
the triangular inequality). In that sense, the (data-driven) symmetric
distance associated with the Oja and Paindaveine \cite
{Oja:2005p311} \emph
{lift-interdirections}, that was recently used to build
nearest-neighbor regression estimators in Biau \textit{et~al.}
\cite{Biau:2012uh}, is
unnecessarily strong. Also, only an ordering of the ``distances'' is
needed to identify nearest neighbors. This \emph{ordering} of
distances \emph{from a fixed point} ${\mathbf{x}}$ is exactly
what the depth-based ${\mathbf{x}}$-outward ordering above is
providing.



\section{Depth-based $k$NN classifiers}
\label{classisec}
In this section, we first define the proposed depth-based classifiers
and present some of their basic properties (Section~\ref
{definprocsec}). We then state the main result of this paper, related
to their consistency  (Section~\ref{consistencysec}).

\subsection{Definition and basic properties}
\label{definprocsec}

The standard $k$-nearest-neighbor ($k$NN) procedure classifies the
point ${\mathbf{x}}$ into Population 1 iff there are more
observations from Population 1 than from Population 0 in the smallest
Euclidean ball centered at ${\mathbf{x}}$ that contains $k$
data points. Depth-based $k$NN classifiers are naturally obtained by
replacing these Euclidean neighborhoods with the depth-based
neighborhoods introduced above, that is, the proposed $k$NN procedure
classifies ${\mathbf{x}}$ into Population 1 iff there are more
observations from Population 1 than from Population 0 in the smallest
depth-based neighborhood of ${\mathbf{x}}$ that contains $k$
observations -- that is, in $R_{{\mathbf{x}}}^{\beta(n)}$,
$\beta=k/n$. In other words, the proposed depth-based classifier is
defined as
%
\begin{equation}
\label{dbclass} \hat{m}_D^{(n)}({\mathbf{x}}) = 
\mathbb{I} \Biggl[ \sum
_{i=1}^n \mathbb{I}[Y_i=1]
W_i^{\beta(n)}({\mathbf{x}}) > \sum
_{i=1}^n \mathbb{I}[Y_i=0]
W_i^{\beta(n)}({\mathbf{x}}) \Biggr],
\end{equation}
with
$
W_i^{\beta(n)}({\mathbf{x}})
=
\frac{1}{K_{\mathbf{x}}^{\beta(n)}}
\mathbb{I}[\Xb_i\in R^{\beta(n)}_{{\mathbf{x}}}]
$,
where $K_{\mathbf{x}}^{\beta(n)}=\sum_{j=1}^n \mathbb
{I}[\Xb_j\in R_{{\mathbf{x}}}^{\beta(n)}]$ still denotes the
number of observations in the depth-based neighborhood $R_{
{\mathbf{x}}}^{\beta(n)}$. Since
%
\begin{equation}
\label{dbclassplugin} \hat{m}_D^{(n)}({\mathbf{x}}) = 
\mathbb{I} \bigl[  \hat{\eta}^{(n)}_D({
\mathbf{x}})> 1/2 \bigr],\qquad  \mbox{with } \hat{\eta}^{(n)}_D({
\mathbf{x}})=\sum_{i=1}^n
\mathbb{I}[Y_i=1] W_i^{\beta(n)}({\mathbf{x}}),
\end{equation}
%
the proposed classifier is actually the one obtained by plugging, in
(\ref{Bayes}), the depth-based estimator $\hat{\eta}^{(n)}
_D({\mathbf{x}})$ of the conditional expectation $\eta
({\mathbf{x}})$. This will be used in the proof of Theorem~\ref{mainresult} below. Note that in the univariate case $(d=1)$,
$\hat m^{(n)}_D$, irrespective of the statistical depth function $D$,
reduces to the standard (Euclidean) $k$NN classifier.

It directly follows from property (P1) that the proposed classifier is
affine-invariant, in the sense that the outcome of the classification
will not be affected if $\Xb_1,\ldots,\Xb_n$ and ${\mathbf{x}}$ are
subject to a common (arbitrary) affine transformation. This
clearly improves over the standard $k$NN procedure that, for example,
is sensitive to unit changes. Of course, one natural way to define an
affine-invariant $k$NN classifier is to apply the original $k$NN
procedure on the standardized data points $\hat{\Sigmab}^{-1/2}\Xb
_i$, $i=1,\ldots,n$, where $\hat{\Sigmab}$ is an affine-equivariant
estimator of shape -- in the sense that
\[
\hat{\Sigmab}(\Ab\Xb_1+\bb,\ldots, \Ab\Xb_n+\bb)
\propto \Ab\hat{\Sigmab}(\Xb_1,\ldots, \Xb_n)
\Ab'
\]
for any invertible $d\times d$ matrix $\Ab$ and any $d$-vector $\bb$.
A natural choice for $\hat{\Sigmab}$ is the regular covariance
matrix, but more robust choices, such as, for example, the shape
estimators from Tyler \cite{Tyler:1987p55},
D{\"u}mbgen \cite{Dumbgen:1998ku}, or Hettmansperger and Randles \cite{Hettmansperger:2002p97} would allow to
get rid of any moment
assumption. Here, we stress that, unlike our \emph{adaptive}
depth-based methodology, such a transformation approach leads to
neighborhoods that do not exploit the geometry of the distribution in
the vicinity of the point ${\mathbf{x}}$ to be classified
(these neighborhoods indeed all are ellipsoids with ${\mathbf
{x}}$-independent orientation and shape); as we show through simulations
below, this results into significantly worse performances.

Most depth-based classifiers available -- among which those relying
on the max-depth approach of Liu, Parelius and Singh \cite
{Liu:1999wm} and Ghosh and Chaudhuri \cite{Ghosh:2005wx}, as
well as the more efficient ones from Li, Cuesta-Albertos and Liu \cite
{Lietal2012} -- suffer from the ``outsider problem\footnote{The term
``outsider'' was recently introduced in Lange, Mosler and Mozharovskyi \cite
{Lange12}.}'': if the
point ${\mathbf{x}}$ to be classified does not sit in the
convex hull of any of the two populations, then most statistical depth
functions will give ${\mathbf{x}}$ zero depth with respect to
each population, so that ${\mathbf{x}}$ cannot be classified
through depth. This is of course undesirable, all the more so that such
a point ${\mathbf{x}}$ may very well be easy to classify. To
improve on this, Hoberg and Mosler \cite{Hoberg2006} proposed
extending the original
depth fields by using the Mahalanobis depth outside the supports of
both populations, a solution that quite unnaturally requires combining
two depth functions. Quite interestingly, our symmetrization
construction implies that the depth-based $k$NN classifier (that
involves one depth function only) does not suffer from the outsider
problem; this is an important advantage over competing depth-based classifiers.

While our depth-based classifiers in (\ref{dbclass}) are perfectly
well-defined and enjoy, as we will show in Section~\ref
{consistencysec} below, excellent consistency properties, practitioners
might find quite arbitrary that a point ${\mathbf{x}}$ such
that ${\sum_{i=1}^n} \mathbb{I}[Y_i=1] W_i^{\beta(n)}(
{\mathbf{x}}) =\sum_{i=1}^n \mathbb{I}[Y_i=0] W_i^{\beta
(n)}({\mathbf{x}})$ is assigned to Population~0. Parallel to
the standard $k$NN classifier, the classification may alternatively be
based on the population of the next neighbor. Since ties are likely to
occur when using depth, it is natural to rather base classification on
the proportion of data points from each population in the next depth
region. Of course, if the next depth region still leads to an ex-aequo,
the outcome of the classification is to be determined on the subsequent
depth regions, until a decision is reached (in the unlikely case that
an ex-aequo occurs for all depth regions to be considered,
classification should then be done by flipping a coin). This treatment
of ties is used whenever real or simulated data are considered below.

Finally, practitioners have to choose some value for the smoothing
parameter $k_n$. This may be done, for example, through
cross-validation (as we will do in the real data example of
Section~\ref{realexsec}).
The value of $k_n$ is likely to have a strong impact on finite-sample
performances, as confirmed in the simulations we conduct in
Section~\ref{simusec}.

\subsection{Consistency results}
\label{consistencysec}

As expected, the local (nearest-neighbor) nature of the proposed
classifiers makes them consistent under very mild conditions. This,
however, requires that the statistical depth function $D$ satisfies the
following further properties:
\begin{enumerate}[(Q3)]
\item[(Q1)] \emph{continuity}: if $P$ is symmetric about $\thetab$
and admits a density that is positive at $\thetab$,
then ${\mathbf{x}}\mapsto D({\mathbf{x}},P)$ is
continuous in a neighborhood of $\thetab$;
\item[(Q2)] \emph{unique maximization at the symmetry center}:
if $P$ is symmetric about $\thetab$ and admits a density that is
positive at $\thetab$,
then $D(\thetab,P)>D({\mathbf{x}},P)$ for all $
{\mathbf{x}}\neq\thetab$;
\item[(Q3)] \emph{consistency}: for any bounded $d$-dimensional Borel
set $B$, $\sup_{{\mathbf{x}}\in B}|D({\mathbf
{x}},P^{(n)})-\linebreak[4] D({\mathbf{x}},P)|=\mathrm{o}(1)$ almost surely as $n\to
\infty$, where $P^{(n)}$ denotes the empirical distribution associated
with $n$ random vectors that are i.i.d. $P$.
\end{enumerate}

Property (Q2) complements property (P2), and, in view of property (P3),
only further requires that $\thetab$ is a strict local maximizer of
${\mathbf{x}}\mapsto D({\mathbf{x}},P)$. Note that
properties (Q1)--(Q2) jointly ensure that the depth-based neighborhoods
of ${\mathbf{x}}$ from Section~\ref{neighsec} collapse to the
singleton $\{{\mathbf{x}}\}$ when the depth level increases to
its maximal value. Finally, since our goal is to prove that our
classifier satisfies an asymptotic property (namely, consistency), it
is not surprising that we need to control the asymptotic behavior of
the sample depth itself (property (Q3)). As shown by Theorem~\ref{Q}, properties (Q1)--(Q3) are satisfied for many classical
depth functions.

We can now state the main result of the paper, that shows that, unlike
their depth-based competitors (that at best are consistent under
semiparametric -- typically elliptical -- distributional assumptions),
the proposed classifiers achieve consistency under virtually any
absolutely continuous distributions. We speak of \emph{nonparametric
consistency}, in order to stress the difference with the stronger \emph
{universal consistency} of the standard $k$NN classifiers.

\begin{Theor}
\label{mainresult}
Let $D$ be a depth function satisfying \textup{(P2)}, \textup{(P3)} and \textup{(Q1)}--\textup{(Q3)}.
Let $k_n$ be a sequence of positive integers such that $k_n\to\infty$
and $k_n=\mathrm{o}(n)$ as $n\to\infty$.
Assume that, for $j=0,1$, $\Xb\mid [Y=j]$ admits a density $f_j$ whose
collection of discontinuity points
has Lebesgue measure zero. Then the depth-based $k_n$NN classifier
$m_D^{(n)}$ in (\ref{dbclass}) is consistent in the sense that
\[
P\bigl[m_D^{(n)}(\Xb)\neq Y \mid \mathcal{D}_n
\bigr]-P\bigl[m_\mathrm{Bayes}(\Xb)\neq Y\bigr]=\mathrm{o}_P(1) \qquad \mbox{as } \ny,
\]
where $\mathcal{D}_n$ is the sigma-algebra associated with $(\Xb
_i,Y_i)$, $i=1,\ldots,n$.
\end{Theor}

Classically, consistency results for classification are based on a
famous theorem from Stone \cite{Stone:1977tk};
see, for example, Theorem~6.3
in Devroye, Gy{\"o}rfi and Lugosi \cite{Devroye:1996to}. However, it is
an open question whether
condition (i) of this theorem holds or not for the proposed
classifiers, at least for some particular statistical depth functions.
A sufficient condition for condition (i) is actually that there exists
a partition of $\R^d$ into cones $C_1,\ldots, C_{\gamma_d}$ with
vertex at the origin of $\R^d$ ($\gamma_d$ not depending on $n$) such
that, for any $\Xb_i$ and any $j$, there exist (with probability one)
at most $k$ data points $\Xb_\ell\in\Xb_i+C_j$ that have $\Xb_i$
among their $k$ depth-based nearest neighbors. Would this be
established for some statistical depth function $D$, it would prove
that the corresponding depth-based $k_n$NN classifier $\hat{m}_D^{(n)}$
is \emph{universally consistent}, in the sense that consistency holds
without \emph{any} assumption on the distribution of $(\Xb,Y)$.

Now, it is clear from the proof of Stone's theorem that this condition
(i) may be dropped if one further assumes that $\Xb$ admits a
uniformly continuous density. This is however a high price to pay, and
that is the reason why the proof of Theorem~\ref{mainresult} rather
relies on an argument recently used in Biau \textit{et~al.} \cite
{Biau:2012uh}; see the \hyperref[appsec]{Appendix}.


\section{Simulations}
\label{simusec}

We performed simulations in order to evaluate the finite-sample
performances of the proposed depth-based $k$NN classifiers. We
considered six setups, focusing on bivariate $\Xb_i$'s ($d=2$) with
equal a priori probabilities ($\pi_0=\pi_1=1/2$), and involving the
following densities $f_0$ and $f_1$:

\begin{Setup}[(Multinormality)] \label{set1}
$f_j$, $j=0,1$, is the pdf of the
bivariate normal distribution with mean vector $\mub_j$ and covariance
matrix $\Sigmab_j$, where
\[
\mub_0= \left( %
\begin{array} {c} 0
\\
0 \end{array} %
\right) , \qquad \mub_1= \left( %
\begin{array}
{c} 1
\\
1 \end{array} %
\right) ,\qquad  \Sigmab_0= \left( %
\begin{array} {c@{\quad}c} 1 & 1
\\
1 & 4 \end{array} %
\right) , \qquad \Sigmab_1=4
\Sigmab_0.
\]
\end{Setup}

\begin{Setup}[(Bivariate Cauchy)]\label{set2}
 $f_j$, $j=0,1$, is the pdf of the
bivariate Cauchy distribution with location center $\mub_j$ and
scatter matrix $\Sigmab_j$, with the same values of $\mub_j$ and
$\Sigmab_j$ as in Setup 1.
\end{Setup}

\begin{Setup}[(Flat covariance structure)] \label{set3}
$f_j$, $j=0,1$, is the
pdf of the bivariate normal distribution with mean vector $\mub_j$ and
covariance matrix $\Sigmab_j$, where
\[
\mub_0= \left( %
\begin{array} {c} 0
\\
0 \end{array} %
\right) , \qquad \mub_1= \left( %
\begin{array} {c} 1
\\
1 \end{array} %
\right) , \qquad \Sigmab_0= \left(
\begin{array} {c@{\quad}c} 5^2 & 0
\\
0 & 1 \end{array} %
\right) , \qquad \Sigmab_1=
\Sigmab_0.
\]
\end{Setup}

\begin{Setup}[(Uniform distributions on half-moons)]\label{set4}
 $f_0$ and
$f_1$ are the densities of
\[
\left( %
\begin{array} {c} X
\\
Y \end{array} %
\right) = \left( %
\begin{array} {c} U
\\
V \end{array} %
\right) \quad \mbox{and}\quad  \left( %
\begin{array}
{c} X
\\
Y \end{array} %
\right) = \left( %
\begin{array} {c} -0.5
\\
2 \end{array} %
\right) + \left( %
\begin{array} {c@{\quad}c} 1 & 0.5
\\
0.5 & -1 \end{array} %
\right) \left( %
\begin{array} {c} U
\\
V \end{array} %
\right),
\]
respectively, where $U\sim\Unif(-1,1)$ and $V|[U=u] \sim
\Unif(1-u^2, 2(1-u^2))$;
\end{Setup}

\begin{Setup}[(Uniform distributions on rings)]\label{set5} $f_0$ and $f_1$ are
the uniform distributions on the concentric rings $\{{\mathbf{x}}\in
\R^2: 1\leq\|{\mathbf{x}}\| \leq2\}$ and $\{
{\mathbf{x}}\in\R^2: 1.75\leq\|{\mathbf{x}}\| \leq
2.5\}$, respectively.
\end{Setup}

\begin{Setup}[(Bimodal populations)]\label{set6}
 $f_j$, $j=0,1$, is the pdf of
the multinormal mixture $\frac{1}{2}\mathcal{N}(\mub^I_{j},\allowbreak \Sigmab
^I_{j})+\frac{1}{2}\mathcal{N}(\mub^{\mathit{II}}_{j},\Sigmab^{\mathit{II}}_{j})$, where
\begin{eqnarray*}
\mub^I_{0}&=& \left( %
\begin{array} {c} 0
\\
0 \end{array} %
\right) , \qquad \mub^{\mathit{II}}_{0}= \left(
\begin{array} {c} 3
\\
3 \end{array} %
\right) , \qquad \Sigmab^I_{0}= \left(
\begin{array} {c@{\quad}c} 1 & 1
\\
1 & 4 \end{array} %
\right) , \qquad \Sigmab^{\mathit{II}}_{0}=4
\Sigmab^I_{0},
\\
\mub^I_{1}&= &\left( %
\begin{array} {c} 1.5
\\
1.5 \end{array} %
\right) ,\qquad  \mub^{\mathit{II}}_{1}= \left(
\begin{array} {c} 4.5
\\
4.5 \end{array} %
\right) , \qquad \Sigmab^I_{1}=
\left( %
\begin{array} {c@{\quad}c} 4 & 0
\\
0 & 0.5 \end{array} %
\right)  \quad \mbox{and}\quad  \Sigmab^{\mathit{II}}_{1}=
\left( %
\begin{array} {c@{\quad}c} 0.75 & 0
\\
0 & 5 \end{array} %
\right) .
\end{eqnarray*}
\end{Setup}

For each of these six setups, we generated 250 training and test
samples of size $n=n_\mathrm{train}=200$ and $n_\mathrm{test}=100$,
respectively, and evaluated the misclassification frequencies of the
following classifiers:
\begin{enumerate}[3.]
\item[1.] The usual LDA and QDA classifiers (LDA/QDA);
\item[2.] The standard Euclidean $k$NN classifiers ($k$NN), with $\beta
=k/n=0.01$, $0.05$, $0.10$ and $0.40$, and the corresponding
``Mahalanobis'' $k$NN classifiers ($k$NNaff) obtained by performing the
Euclidean $k$NN classifiers on standardized data, where standardization
is based on the regular covariance matrix estimate of the whole
training sample;
\item[3.] The proposed depth-based $k$NN classifiers (D-$k$NN) for each
combination of the $k$ used in $k$NN$/$$k$NNaff and a statistical depth
function (we focused on halfspace depth, simplicial depth, and
Mahalanobis depth);
\item[4.] The depth vs depth (DD) classifiers from Li, Cuesta-Albertos and Liu
\cite{Lietal2012}, for
each combination of a polynomial curve of degree $m$ ($m=1$, 2, or 3)
and a statistical depth function (halfspace depth, simplicial depth, or
Mahalanobis depth). Exact DD-classifiers (DD) as well as smoothed
versions (DDsm) were actually implemented  -- although, for
computational reasons, only the smoothed version was considered for
$m=3$. Exact classifiers search for the best separating polynomial
curve $(d,r(d))$ of order $m$ passing through the origin and $m$
``DD-points'' $(D_0^{(n)}(\Xb_i),D_1^{(n)}(\Xb_i))$ (see the \hyperref[introsec]{Introduction})
in the sense that it minimizes the misclassification error
%
\begin{equation}
\label{missclass} \sum_{i=1}^n \bigl(
\mathbb{I}[Y_i=1]\mathbb{I}\bigl[ d_i^{(n)}>0
\bigr]+\mathbb {I}[Y_i=0]\mathbb{I}\bigl[-d_i^{(n)}>0
\bigr] \bigr),
\end{equation}
with $d_i^{(n)}:=r(D_0^{(n)}(\Xb_i))-D_1^{(n)}(\Xb_i)$. Smoothed
versions use
derivative-based methods to find a polynomial minimizing (\ref
{missclass}), where the indicator $\mathbb{I}[d>0]$ is replaced by the
logistic function $1/(1+\mathrm{e}^{-td})$ for a suitable $t$. As suggested in
Li, Cuesta-Albertos and Liu \cite{Lietal2012}, value $t=100$ was chosen
in these simulations.
$100$ randomly chosen polynomials were used as starting points for the
minimization algorithm, the classifier using the resulting polynomial
with minimal misclassification (note that this time-consuming scheme
always results into better performances than the one adopted in
Li, Cuesta-Albertos and Liu \cite{Lietal2012}, where only one
minimization is performed, starting from
the best random polynomial considered).
\end{enumerate}

Since the DD classification procedure is a refinement of the max-depth
procedures of Ghosh and Chaudhuri \cite{Ghosh:2005wx} that
leads to better
misclassification rates (see Li, Cuesta-Albertos and Liu \cite
{Lietal2012}), the original
max-depth procedures were omitted in this study.
%
\begin{figure}

\includegraphics{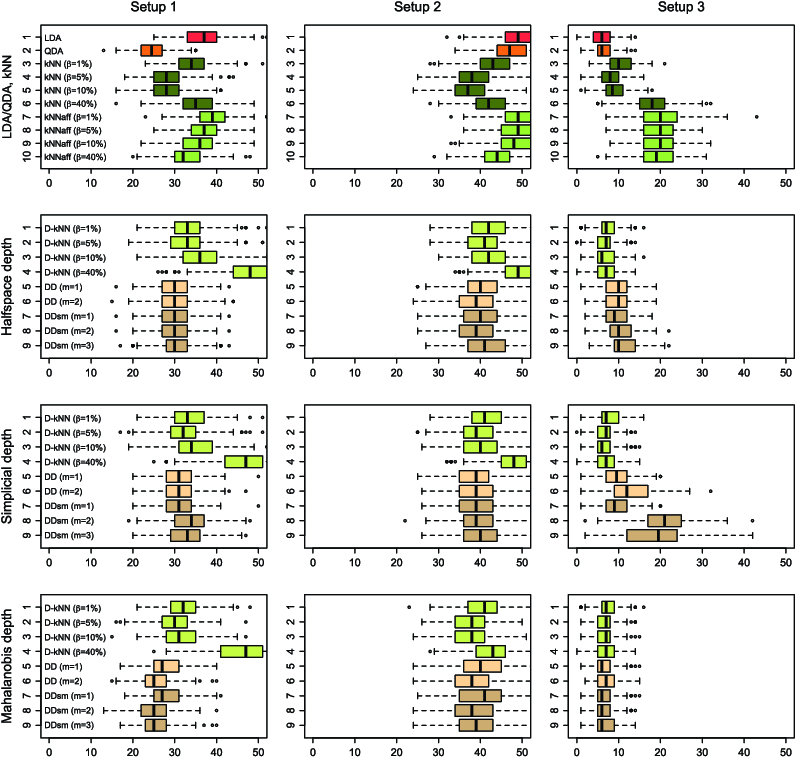}

\caption{Boxplots of misclassification frequencies (in percentages),
from 250 %
replications of Setups \protect\ref{set1} to \protect\ref{set3} described in Section \protect\ref
{simusec}, with
training sample size $n=n_\mathrm{train}=200$ and test sample size
$n_\mathrm{test}=100$, of %
the LDA/QDA classifiers, the Euclidean $k$NN classifiers ($k$NN) and
their Mahalanobis (affine-invariant) counterparts ($k$NNaff), the
proposed depth-based $k$NN %
classifiers (D-$k$NN), and some exact and smoothed version of the
DD-classifiers (DD and DDsm); see Section \protect\ref{simusec} for details.}
\label{Fig1}
\end{figure}
%
%
\begin{figure}

\includegraphics{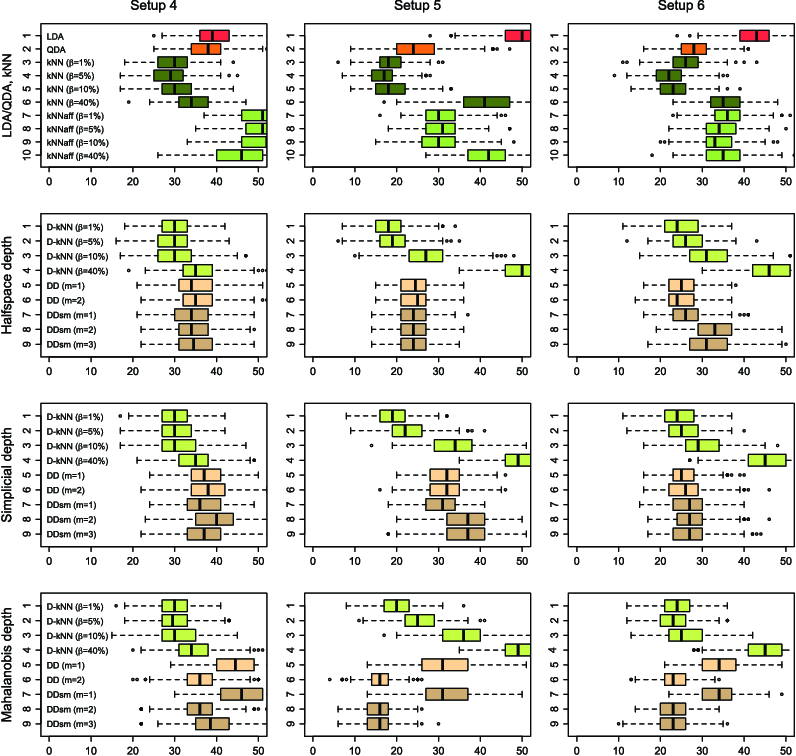}

\caption{Boxplots of misclassification frequencies (in percentages),
from 250
replications of Setups \protect\ref{set4} to \protect\ref{set6} described in Section \protect\ref
{simusec}, with
training sample size $n=n_\mathrm{train}=200$ and test sample size
$n_\mathrm{test}=100$, of %
the LDA/QDA classifiers, the Euclidean $k$NN classifiers ($k$NN) and
their Mahalanobis (affine-invariant) counterparts ($k$NNaff), the
proposed depth-based $k$NN %
classifiers (D-$k$NN), and some exact and smoothed version of the
DD-classifiers (DD and DDsm); see Section \protect\ref{simusec} for details.}
\label{Fig2}
\end{figure}

Boxplots of misclassification frequencies (in percentages) are reported
in Figures~\ref{Fig1} and~\ref{Fig2}. 
It is seen that
in most setups, the proposed depth-based $k$NN classifiers compete well
with the Euclidean $k$NN classifiers. The latter, however, should be
avoided since (i) their outcome may unpleasantly depend on measurement
units, and since (ii) the spherical nature of the neighborhoods used
lead to performances that are severely affected by the -- notoriously
delicate -- choice of $k$; see the ``flat'' Setup \ref{set3}.
This motivates restricting to affine-invariant classifiers, that (i)
are totally insensitive to any unit changes and that (ii) can adapt to
the flat structure of Setup \ref{set3} as they show there performances that are
much more stable in $k$.

Now, regarding the comparisons between affine-invariant classifiers,
the simulations results lead to the following conclusions:
(i)
the proposed affine-invariant depth-based classifiers outperform the
natural affine-invariant versions of $k$NN classifiers. In other words,
the natural way to make the standard $k$NN classifier affine-invariant
results into a dramatic cost in terms of finite-sample performances.
%
%
(ii) The proposed depth-based $k$NN classifiers also compete well with
DD-classifiers both in elliptical and non-elliptical setups.
Away from ellipticity (Setups \ref{set4} to \ref{set6}), in particular, they perform at
least as well -- and sometimes outperform (Setup \ref{set4}) -- DD-classifiers; a
single exception is associated with the use of Mahalanobis depth in
Setup \ref{set5}, where the DD-classifiers based on $m=2,3$ perform better.
Apparently, another advantage of depth-based $k$NN classifiers over
DD-classifiers is that their finite-sample performances depend much
less on the statistical depth function $D$ used.
%


%
\section{Real-data examples}
\label{realexsec}

In this section, we investigate the performances of our depth-based
$k$NN classifiers on two well known benchmark datasets. The first
example is taken from 
Ripley \cite{Ripley:1996} and can be found on the
book's website (\url{http://www.stats.ox.ac.uk/pub/PRNN}). This data set involves
well-specified training and test samples, and we therefore simply
report the test set misclassification rates of the different
classifiers included in the study.
The second example, blood transfusion data, is available at \url{http://archive.ics.uci.edu/ml/index.html}. Unlike the first data set,
no clear partition into a training sample and a test sample is provided
here. As suggested in Li, Cuesta-Albertos and Liu \cite{Lietal2012}, we
randomly performed such a
partition $100$ times (see the details below) and computed the average
test set misclassification rates, together with standard deviations.

A brief description of each dataset is as follows:

\textit{Synthetic data} was introduced and studied in Ripley \cite{Ripley:1996}. The dataset is made of observations
from two
populations, each of them being actually a mixture of two bivariate
normal distributions differing only in location. As mentioned above, a
partition into a training sample and a test sample is provided: the
training and test samples contain $250$ and $1000$ observations,
respectively, and both samples are divided equally between the two populations.
%
\begin{table}[b]
\tablewidth=310pt
\tabcolsep=0pt
\caption{Misclassification rates (for synthetic data) and sample averages and standard deviations
(in parentheses) of misclassification rates obtained from 100 random partitions of the data into
training and test samples (for transfusion data)}\label{Tablereal}
%
\begin{tabular*}{310pt}{@{\extracolsep{\fill}}lll@{}}
\hline
&Synthetic & Transfusion\\
\hline
LDA&10.8& 29.60 (0.9)\\
QDA&10.2& 29.21 (1.5)\\
$k$NN&\phantom{1}{8.7}&29.74 (2.0)\\
$k$NNaff&11.7&30.11 (2.1)\\
D$_H$-$k$NN &10.1&27.75 (1.6)\\
D$_M$-$k$NN &14.4&27.36 (1.5)\\
DD$_H$ ($m=1$)&13.4 &28.26 (1.7)\\
DD$_H$ ($m=2$)&12.9 &28.33 (1.6)\\
DD$_M$ ($m=1$)&17.5&31.44 (0.1)\\
DD$_M$ ($m=2$)&12.0&31.54 (0.6)\\
\hline
\end{tabular*}
\end{table}

\textit{Transfusion data} contains the information on 748 blood donors
selected from the blood donor database of the Blood Transfusion Service
Center in Hsin-Chu City, Taiwan. It was studied in Yeh, Yang and Ting \cite{Yeh:2009}.
The classification problem at hand is to know whether or not the donor
gave blood in March 2007. In this dataset, prior probabilities are not
equal; out of 748 donors, 178 gave blood in March 2007, when 570 did
not. Following Li, Cuesta-Albertos and Liu \cite{Lietal2012}, one out of
two linearly correlated
variables was removed and three measurements were available for each
donor: Recency (number of months since the last donation), Frequency
(total number of donations) and Time (time since the first donation).
The training set consists in $100$ donors from the first class and
$400$ donors from the second, while the rest is assigned to the test
sample (therefore containing $248$ individuals).

Table~\ref{Tablereal} reports the -- exact (synthetic) or averaged
(transfusion) -- misclassification rates of the following classifiers:
the linear (LDA) and quadratic (QDA) discriminant rules,
the standard $k$NN classifier ($k$NN) and its Mahalanobis
affine-invariant version ($k$NNaff),
the depth-based $k$NN classifiers using halfspace depth (D$_H$-$k$NN)
and Mahalanobis depth (D$_M$-$k$NN),
and the exact DD-classifiers for any combination of a polynomial order
$m\in\{1,2\}$ and a statistical depth function among the two
considered for depth-based $k$NN classifiers, namely the halfspace
depth (DD$_H$) and the Mahalanobis depth (DD$_M$) -- smoothed
DD-classifiers were excluded from this study, as their performances,
which can only be worse than those of exact versions, showed much
sensitivity to the smoothing parameter $t$; see Section~\ref{simusec}.
For all nearest-neighbor classifiers, leave-one-out cross-validation
was used to determine $k$.

The results from Table~\ref{Tablereal} indicate that depth-based $k$NN
classifiers perform very well in both examples. For synthetic data, the
halfspace depth-based $k$NN classifier (10.1\%) is only dominated by
the standard (Euclidean) $k$NN procedure (8.7\%). The latter, however,
has to be discarded as it is dependent on scale and shape changes -- in
line with this, note that the ``$k$NN classifier'' applied in
Dutta and Ghosh \cite{Dutta:wt} is actually the $k$NNaff
classifier (11.7\%), as
classification in that paper is performed on standardized data. The
Mahalanobis depth-based $k$NN classifiers (14.4\%) does not perform as
well as its halfspace counterpart. For transfusion data, however, both
depth-based $k$NN classifiers dominate their competitors.



\section{Final comments}
\label{finalsec}

The depth-based neighborhoods we introduced are of interest in other
inference problems as well. As an illustration, consider the regression
problem where the conditional mean function
$
{\mathbf{x}}\mapsto m({\mathbf{x}})=E[Y \mid \Xb
={\mathbf{x}}]
$
is to be estimated on the basis of mutually independent copies $(\Xb
_i,Y_i)$, $i=1,\ldots,n$ of a random vector $(\Xb,Y)$ with values in
$\R^d\times\R$, or the problem of estimating the common density $f$
of \mbox{i.i.d.} random $d$-vectors $\Xb_i$, $i=1,\ldots,n$. The
classical $k_n$NN estimators for these problems are
%
\begin{eqnarray}
\label{denskNN}  
\hat{m}^{(n)}({\mathbf{x}}) &=& \sum_{i=1}^n
W_i^{\beta_n(n)}({\mathbf{x}}) Y_i =
\frac{1}{k_n} \sum_{i=1}^n \mathbb{I}
\bigl[\Xb_i\in B_{\mathbf{x}}^{\beta_n{(n)}} \bigr]
Y_i,\quad \mbox{and}
\nonumber
\\[-8pt]
\\[-8pt]
\nonumber
\hat f^{(n)}({\mathbf{x}}) &=&\frac{k_n}{n \mu_d(B_{\mathbf{x}}^{\beta_n{(n)}})}
\end{eqnarray}
where $\beta_n=k_n/n$, $B_{\mathbf{x}}^{\beta(n)}$ is the
smallest Euclidean ball centered at ${\mathbf{x}}$ that
contains a proportion $\beta$ of the $\Xb_i$'s, and $\mu_d$ stands
for the Lebesgue measure on $\R^d$. Our construction naturally leads
to considering the depth-based $k_n$NN estimators $\hat m^{(n)}
_D({\mathbf{x}})$ and $\hat{f}^{(n)}_D({\mathbf{x}})$
obtained by replacing in (\ref{denskNN}) the Euclidean neighborhoods
$B_{\mathbf{x}}^{\beta_n}$ with their depth-based
counterparts $R_{\mathbf{x}}^{\beta_n{(n)}}$ and $k_n=\sum_{i=1}^n \mathbb{I} [\Xb_i\in B_{\mathbf{x}}^{\beta_n{(n)}}
]$ with $K_{\mathbf{x}}^{\beta_n(n)}=\sum_{i=1}^n
\mathbb{I} [\Xb_i\in R_{\mathbf{x}}^{\beta_n{(n)}}  ]$.

A thorough investigation of the properties of these depth-based
procedures is of course beyond the scope of the present paper. It is,
however, extremely likely that the excellent consistency properties
obtained in the classification problem extend to these nonparametric
regression and density estimation setups. Now, recent works in density
estimation indicate that using non-spherical (actually, ellipsoidal)
neighborhoods may lead to better finite-sample properties; see, for
example, Chac{\'o}n \cite{Chacon:2009cq} or
Chac{\'o}n, Duong and Wand \cite{Chacon:2011}. In that respect,
the depth-based $k$NN estimators above are very promising since they
involve non-spherical (and for most classical depth, even
non-ellipsoidal) neighborhoods whose shape is determined by the local
geometry of the sample. Note also that depth-based neighborhoods only
require choosing a single scalar bandwidth parameter (namely, $k_n$),
whereas general $d$-dimensional ellipsoidal neighborhoods impose
selecting $d(d+1)/2$ bandwidth parameters.


\begin{appendix}

\section*{Appendix: Proofs}
\label{appsec}

The main goal of this appendix is to prove Theorem~\ref{mainresult}.
We will need the following lemmas.

\begin{Lem}
\label{lemme1}
Assume that the depth function $D$ satisfies \textup{(P2)}, \textup{(P3)}, \textup{(Q1)}, and \textup{(Q2)}.
Let $P$ be a probability measure that is symmetric about $\thetab$ and
admits a density that is positive at $\thetab$.
Then,
\textup{(i)} for all $a>0$, there exists $\alpha<\alpha_*=\max_{
{\mathbf{x}}\in\R^d}D({\mathbf{x}},P)$ such that $R_\alpha
(P)\subset B_{\thetab}(a):=\{{\mathbf{x}}\in\R^d:\|
{\mathbf{x}}-\thetab\|\leq a\}$;
\textup{(ii)} for all $\alpha<\alpha_*$, there exists $\xi>0$ such that
$B_{\thetab}(\xi)\subset R_\alpha(P)$.
\end{Lem}

\begin{pf}
(i) First, note that the existence of $\alpha_*$ follows from property
(P2). Fix then $\delta>0$ such that ${\mathbf{x}}\mapsto
D({\mathbf{x}},P)$ is continuous over $B_{\thetab}(\delta)$;
existence of $\delta$ is guaranteed by property (Q1). Continuity
implies that ${\mathbf{x}}\mapsto D({\mathbf{x}},P)$
reaches a minimum in $B_{\thetab}(\delta)$, and property (Q2) entails
that this minimal value, $\alpha_\delta$ say, is strictly smaller
than $\alpha_*$. Using property (Q1) again, we obtain that, for each
$\alpha\in[\alpha_\delta,\alpha_*]$,
\begin{eqnarray*}
r_\alpha\dvtx   \mathcal{S}^{d-1} & \to&\R^+,
\\
 \mathbf{u} &\mapsto&\sup\bigl\{ r\in\R^+: \thetab+r \ub\in R_{\alpha}(P)\bigr\}
\end{eqnarray*}
is a continuous function that converges pointwise to $r_{\alpha_*}(\ub
)\equiv0$ as $\alpha\to\alpha_*$. Since $\mathcal{S}^{d-1}$ is
compact, this convergence is actually uniform, that is, $\sup_{\ub\in
\mathcal{S}^{d-1}} | r_\alpha(\ub) |=\mathrm{o}(1)$ as $\alpha\to\alpha
_*$. Part (i) of the result follows.

(ii) Property (Q2) implies that, for any $\alpha\in[\alpha_\delta
,\alpha_*)$, the mapping $r_\alpha$ takes values in $\R_0^+$.
Therefore, there exists $\ub_0(\alpha)\in\mathcal{S}^{d-1}$ such
that $r_\alpha(\ub)\geq r_\alpha(\ub_0(\alpha))=\xi_\alpha>0$.
This implies that, for all $\alpha\in[\alpha_\delta,\alpha_*)$,
we have $B_{\thetab}(\xi_\alpha)\subset R_\alpha(P)$, which proves
the result for these values of $\alpha$. Nestedness of the $R_\alpha
(P)$'s, which follows from property (P3), then establishes the result
for an arbitrary $\alpha<\alpha_*$.
\end{pf}

\begin{Lem}
\label{lemme2}
Assume that the depth function $D$ satisfies \textup{(P2)}, \textup{(P3)}, and \textup{(Q1)}--\textup{(Q3)}.
Let $P$ be a probability measure that is symmetric about $\thetab$ and
admits a density that is positive at $\thetab$.
%
Let $\Xb_1,\ldots,\Xb_n$ be \mbox{i.i.d.} $P$ and denote by $\Xb
_{\thetab,(i)}$ the $i$th depth-based nearest neighbor of $\thetab$.
Let $K_{\thetab}^{\beta_n(n)}$ be the number of depth-based nearest
neighbors in $R_{\thetab}^{\beta_n}(P^{(n)})$, where $\beta_n=k_n/n$ is
based on a sequence $k_n$ that is as in Theorem~\ref{mainresult} and
$P^{(n)}$ stands for the empirical distribution of $\Xb_1,\ldots,\Xb_n$.
Then,
for any $a>0$, there exists $n=n(a)$ such that $\sum_{i=1}^{K^{\beta
_n(n)}_{\thetab}} \mathbb{I}[\|\Xb_{\thetab,(i)}-\thetab\|>a]=0$
almost surely for all $n\geq n(a)$.
\end{Lem}

Note that, while $\Xb_{\thetab,(i)}$ may not be properly defined
(because of ties), the quantity\linebreak[4]  $\sum_{i=1}^{K^{\beta_n(n)}_{\thetab}}
\mathbb{I}[\|\Xb_{\thetab,(i)}-\thetab\|>a]=0$ always is.
\begin{pf*}{Proof of Lemma~\ref{lemme2}}
Fix $a>0$.
By Lemma~\ref{lemme1}, there exists $\alpha<\alpha_*$ such that
$R_\alpha(P)\subset B_{\thetab}(a)$.
Fix then $\bar{\alpha}$ and $\varepsilon>0$ such that $\alpha<\bar
{\alpha}-\varepsilon<\bar{\alpha}+\varepsilon<\alpha_*$. Theorem~4.1 in Zuo and Serfling \cite{Zuo:2000p244} and the fact
that $P^{(n)}_{\thetab}\to
P_{\thetab}=P$ weakly as $n\to\infty$ (where $P^{(n)}_{\thetab}$ and
$P_{\thetab}$ are the $\thetab$-symmetrized versions of $P^{(n)}$ and $P$,
respectively) then entail that there exists an integer $n_0$ such that
\[
R_{\bar{\alpha}+\varepsilon}(P) \subset R_{\bar{\alpha}}\bigl(P^{(n)}_{\thetab}
\bigr) \subset R_{\bar{\alpha}-\varepsilon}(P) \subset R_\alpha(P)
\]
almost surely for all $n\geq n_0$. From Lemma~\ref{lemme1} again,
there exists $\xi>0$ such that $B_{\thetab}(\xi)\subset R_{\bar
{\alpha}+\varepsilon}(P)$. Hence, for any $n\geq n_0$, one has that
\[
B_{\thetab}(\xi) \subset R_{\bar{\alpha}}\bigl(P^{(n)}_{\thetab}
\bigr) \subset B_{\thetab}(a)
\]
almost surely.

Putting $N_n=\sum_{i=1}^n \mathbb{I}[\Xb_i\in B_{\thetab}(\xi)]$,
the SLLN yields that $N_n/n\to P[B_{\thetab}(\xi)]=\break P[B_{\thetab}(\xi
)]>0$ as $\ny$, since $\Xb\sim P$ admits a density that, from
continuity, is positive over a neighborhood of $\thetab$. Since
$k_n=\mathrm{o}(n)$ as $n\to\infty$, this implies that, for all $n\geq\tilde
n_0\ (\geq n_0)$,
\[
\sum_{i=1}^n \mathbb{I}\bigl[
\Xb_i\in R_{\bar{\alpha}}\bigl(P^{(n)}_{\thetab}\bigr)
\bigr] \geq N_n \geq k_n
\]
almost surely. It follows that, for such values of $n$,
\[
R_{\thetab}^{\beta_n}\bigl(P^{(n)}\bigr) =R^{\beta_n}
\bigl(P^{(n)}_{\thetab}\bigr) \subset R_{\bar{\alpha}}
\bigl(P_{\thetab}^{(n)}\bigr) \subset B_{\thetab}(a)
\]
almost surely, with $\beta_n=k_n/n$. Therefore, $\max_{i=1,\ldots
,K^{\beta_n(n)}_{\thetab}}\| \Xb_{\thetab,(i)} - \thetab\|\leq a$
almost surely for large $n$, which yields the result.
\end{pf*}

\begin{Lem}
\label{lemme3}
For a ``plug-in'' classification rule $\tilde{m}^{(n)}({\mathbf
{x}})=\mathbb{I}[\tilde{\eta}^{(n)}({\mathbf{x}})>1/2]$
obtained from a regression estimator $\tilde{\eta}^{(n)}(
{\mathbf{x}})$ of $\eta({\mathbf{x}})
=E[\mathbb{I}[Y=1] \mid \Xb={\mathbf{x}}]$, one has that
$
P[\tilde{m}^{(n)}(\Xb)\neq Y] - L_\mathrm{opt}
\leq
2  ( E[ (\tilde{\eta}^{(n)}(\Xb)-\eta(\Xb))^2 ]  )^{1/2}
$,
where $L_\mathrm{opt}=P[m_\mathrm{Bayes}(\Xb)\neq Y]$ is the
probability of misclassification of the Bayes rule.
\end{Lem}

\begin{pf}
Corollary~6.1 in Devroye, Gy{\"o}rfi and Lugosi \cite{Devroye:1996to}
states that
\[
P\bigl[\tilde{m}^{(n)}(\Xb)\neq Y \mid \mathcal{D}_n\bigr] -
L_\mathrm{opt} \leq 2 E\bigl[ \bigl|\tilde{\eta}^{(n)}(\Xb)-\eta(\Xb)\bigr|
\mid \mathcal{D}_n \bigr],
\]
where $\mathcal{D}_n$ stands for the sigma-algebra associated with the
training sample $(\Xb_i,Y_i)$, $i=1,\ldots,n$. Taking expectations in
both sides of this inequality and applying Jensen's inequality readily
yields the result.
\end{pf}
\begin{pf*}{Proof of Theorem~\ref{mainresult}}
From Bayes' theorem, $\Xb$ admits the density ${\mathbf{x}}\mapsto
f({\mathbf{x}})=\pi_0f_0({\mathbf{x}})+\pi_1f_1({\mathbf{x}})$.
Letting $\Supp_+(f)=\{
{\mathbf{x}}\in\R^d:f({\mathbf{x}})>0\}$ and
writing $C(f_j)$ for the collection of continuity points of $f_j$,
$j=0,1$, put $N=\Supp_+(f)\cap C(f_0)\cap C(f_1)$. Since, by
assumption, $\R^d\setminus C(f_j)$ ($j=0,1$) has Lebesgue measure
zero, we have that
\begin{eqnarray*}
P\bigl[\Xb\in\R^d\setminus N\bigr] &\leq& P\bigl[\Xb\in
\R^d\setminus\Supp_+(f)\bigr] + \sum_{j\in\{0,1\}}P
\bigl[\Xb\in\R^d\setminus C(f_j)\bigr]
\\
&=& \int_{\R^d\setminus\Supp_+(f)} f({\mathbf{x}}) \,\mathrm{d}x =0,
\end{eqnarray*}
so that $P[\Xb\in N]=1$.
Note also that ${\mathbf{x}}\mapsto\eta({\mathbf{x}})=\pi
_1f_1({\mathbf{x}})/(\pi_0f_0({\mathbf{x}})+\pi_1f_1({\mathbf
{x}}))$ is continuous over $N$.

Fix ${\mathbf{x}}\in N$ and let $Y_{{\mathbf{x}},(i)}=Y_{j({\mathbf
{x}})}$ with $j({\mathbf{x}})$
such that $\Xb_{{\mathbf{x}},(i)}=\Xb_{j({\mathbf{x}})}$. With this
notation, the estimator $\hat{\eta}^{(n)}
_D({\mathbf{x}})$ from Section~\ref{definprocsec} rewrites
\[
\hat{\eta}^{(n)}_D({\mathbf{x}}) = \sum
_{i=1}^n Y_i W_i^{\beta(n)}({
\mathbf{x}}) = \frac{1}{K^{\beta(n)}_{\mathbf{x}}} \sum_{i=1}^{K^{\beta
(n)}_{\mathbf{x}}}
Y_{{\mathbf{x}},(i)}.
\]
Proceeding as in Biau \textit{et~al.} \cite{Biau:2012uh}, we
therefore have that (writing
for simplicity $\beta$ instead of $\beta_n$ in the rest of the proof)
\[
T^{(n)}({\mathbf{x}}):=E\bigl[\bigl(\hat{\eta}^{(n)}_D({
\mathbf{x}})-\eta ({\mathbf{x}}) \bigr)^2\bigr] \leq
2T^{(n)}_1({\mathbf{x}})+2T^{(n)}_2({
\mathbf{x}}),
\]
with
\[
T^{(n)}_1({\mathbf{x}}) = E \Biggl[ \Biggl| \frac{1}{K^{\beta(n)}_{\mathbf{x}}}
\sum_{i=1}^{K^{\beta(n)}_{\mathbf{x}}} \bigl( Y_{
{\mathbf{x}},(i)} -
\eta(\Xb_{{\mathbf{x}},(i)}) \bigr) \Biggr|^2 \Biggr]
\]
and
\[
T^{(n)}_2({\mathbf{x}}) = E \Biggl[\Biggl | \frac{1}{K^{\beta(n)}_{\mathbf{x}}} \sum
_{i=1}^{K^{\beta(n)}_{\mathbf{x}}} \bigl( \eta(\Xb _{{\mathbf{x}},(i)})
-\eta({\mathbf{x}}) \bigr) \Biggr|^2 \Biggr].
\]
Writing $\mathcal{D}_X^{(n)}$ for the sigma-algebra generated by $\Xb
_i$, $i=1,\ldots,n$, note that, conditional on $\mathcal{D}_X^{(n)}$,
the $Y_{{\mathbf{x}},(i)} - \eta(\Xb_{{\mathbf{x}},(i)}) $'s,
$i=1,\ldots,n$, are zero mean mutually independent
random variables. Consequently,
\begin{eqnarray*}
T^{(n)}_1({\mathbf{x}}) &=& E \Biggl[ \frac{1}{(K^{\beta(n)}_{\mathbf{x}})^2}
\sum_{i,j=1}^{K^{\beta(n)}_{\mathbf{x}}} E \bigl[ \bigl(
Y_{{\mathbf{x}},(i)} - \eta(\Xb_{{\mathbf{x}},(i)}) \bigr) \bigl( Y_{{\mathbf{x}},(j)} -
\eta(\Xb_{{\mathbf{x}},(j)}) \bigr) \mid \mathcal{D}_X^{(n)} \bigr]
\Biggr]
\\
&=& E \Biggl[ \frac{1}{(K^{\beta(n)}_{\mathbf{x}})^2} \sum_{i=1}^{K^{\beta(n)}_{\mathbf{x}}}
E \bigl[ \bigl( Y_{{\mathbf{x}},(i)} - \eta(\Xb_{{\mathbf{x}},(i)}) \bigr)^2
\mid \mathcal{D}_X^{(n)} \bigr] \Biggr]
\\
&\leq& E \biggl[ \frac{4}{K^{\beta(n)}_{\mathbf{x}}} \biggr] \leq\frac{4}{k_n} =\mathrm{o}(1),
\end{eqnarray*}
as $n\to\infty$, where we used the fact that $K_{\mathbf{x}}^{\beta
(n)}\geq k_n$ almost surely. As for $T^{(n)}_2(
{\mathbf{x}})$, the Cauchy--Schwarz inequality yields (for an arbitrary $a>0$)
\begin{eqnarray*}
T^{(n)}_2({\mathbf{x}}) &\leq& E \Biggl[ \frac{1}{K^{\beta(n)}_{\mathbf{x}}}
\sum_{i=1}^{K^{\beta
(n)}_{\mathbf{x}}} \bigl( \eta(
\Xb_{{\mathbf{x}},(i)}) -\eta({\mathbf{x}}) \bigr)^2 \Biggr]
\\
&=& E \Biggl[ \frac{1}{K^{\beta(n)}_{\mathbf{x}}} \sum_{i=1}^{K^{\beta
(n)}_{\mathbf{x}}}
\bigl( \eta(\Xb_{{\mathbf{x}},(i)}) -\eta({\mathbf{x}}) \bigr)^2
\mathbb{I}\bigl[\|\Xb_{{\mathbf{x}},(i)}-{\mathbf{x}}\|\leq a\bigr] \Biggr]
\\
& &{} + E \Biggl[ \frac{1}{K^{\beta(n)}_{\mathbf{x}}} \sum_{i=1}^{K^{\beta
(n)}_{\mathbf{x}}}
\bigl( \eta(\Xb_{{\mathbf{x}},(i)}) -\eta({\mathbf{x}}) \bigr)^2
\mathbb{I}\bigl[\|\Xb_{{\mathbf{x}},(i)}-{\mathbf{x}}\|>a\bigr] \Biggr]
\\
&\leq& \sup_{{\mathbf{x}}'\in B_{\mathbf{x}}(a)} \bigl|\eta \bigl({\mathbf{x}}'\bigr)
-\eta({\mathbf{x}})\bigr|^2 
%
+ 4 E
\Biggl[ \frac{1}{K^{\beta(n)}_{\mathbf{x}}} \sum_{i=1}^{K^{\beta
(n)}_{\mathbf{x}}}
\mathbb{I}\bigl[\|\Xb_{{\mathbf{x}},(i)}-{\mathbf{x}}\|>a\bigr] \Biggr]
\\
&=:& \tilde T_{2}({\mathbf{x}};a)+\bar T^{(n)}_{2}({
\mathbf{x}};a).
\end{eqnarray*}
Continuity of $\eta$ at ${\mathbf{x}}$ implies that, for any
$\varepsilon>0$, one may choose $a=a(\varepsilon)>0$ so that $\tilde
T_{2}({\mathbf{x}};a(\varepsilon))<\varepsilon$. Since
Lemma~\ref{lemme2} readily yields that $T^{(n)}_{2}({\mathbf
{x}};a(\varepsilon))=0$ for large $n$, we conclude that $T^{(n)}
_2({\mathbf{x}})$ -- hence also $T^{(n)}({\mathbf{x}})$ -- is $\mathrm{o}(1)$.
The Lebesgue dominated convergence theorem then
yields that $E[(\hat{\eta}^{(n)}_D(\Xb)-\eta(\Xb) )^2]
$ is $\mathrm{o}(1)$. Therefore, using the fact that $P[\hat{m}_D^{(n)}(\Xb
)\neq
Y \mid \mathcal{D}_n]\geq L_\mathrm{opt}$ almost surely and applying
Lemma~\ref{lemme3}, we obtain
\begin{eqnarray*}
E \bigl[ \bigl| P\bigl[\hat{m}_D^{(n)}(\Xb)\neq Y \mid
\mathcal{D}_n\bigr]-L_\mathrm{opt}\bigr | \bigr] &=& E\bigl[ P\bigl[
\hat{m}_D^{(n)}(\Xb)\neq Y \mid \mathcal{D}_n
\bigr]-L_\mathrm{opt} \bigr]
\\
&  =& P\bigl[\hat{m}_D^{(n)}(\Xb)\neq Y
\bigr]-L_\mathrm{opt} \leq 2 \bigl( E\bigl[ \bigl(\hat{
\eta}_D^{(n)}(\Xb)-\eta(\Xb)\bigr)^2 \bigr]
\bigr)^{1/2} \\
&=&\mathrm{o}(1),
\end{eqnarray*}
as $\ny$, which establishes the result.
\end{pf*}

Finally, we show that properties (Q1)--(Q3) hold for several classical
statistical depth functions.

\begin{Theor}
\label{Q}
Properties \textup{(Q1)}--\textup{(Q3)} hold for \textup{(i)} the halfspace depth and \textup{(ii)} the
simplicial depth. \textup{(iii)} If the location and scatter functionals $\mub
(P)$ and $\Sigmab(P)$ are such that \textup{(a)} $\mub(P)=\thetab$ as soon as
the probability measure $P$ is symmetric about $\thetab$ and such that
\textup{(b)} the empirical versions $\mub(P^{(n)})$ and $\Sigmab(P^{(n)})$
associated with an i.i.d. sample $\Xb_1,\ldots,\Xb_n$ from
$P$ are strongly consistent for $\mub(P)$ and $\Sigmab(P)$, then
properties \textup{(Q1)}--\textup{(Q3)} also hold for the Mahalanobis depth.
\end{Theor}

\begin{pf}
(i) The continuity of $D$ in property (Q1) actually holds under the
only assumption that $P$ admits a density with respect to the Lebesgue
measure; see Proposition~4 in Rousseeuw and Ruts \cite
{Rousseeuw:1999p158}. Property (Q2)
is a consequence of Theorems 1 and 2 in Rousseeuw and Struyf \cite{Rousseeuw:2004p253} and
the fact that the angular symmetry center is unique for absolutely
continuous distributions; see Serfling \cite
{Serfling:2006up}. For halfspace
depth, property (Q3) follows from (6.2) and (6.6) in Donoho and Gasko \cite{Donoho:1992p214}.

(ii)
The continuity of $D$ in property (Q1) actually holds under the only
assumption that $P$ admits a density with respect to the Lebesgue
measure; see Theorem~2 in Liu \cite{Liu:1990p224}.
Remark C in Liu \cite{Liu:1990p224} shows that,
for an angularly symmetric
probability measure (hence also for a centrally symmetric probability
measure) admitting a density, the symmetry center is the unique point
maximizing simplicial depth provided that the density remains positive
in a neighborhood of the symmetry center; property (Q2) trivially follows.
property (Q3) for simplicial depth is stated in Corollary~1 of
D{\"u}mbgen \cite{Dumbgen:1992tj}.

(iii) This is trivial.
\end{pf}

Finally, note that properties (Q1)--(Q3) also hold for projection depth
under very mild assumptions on the univariate location and scale
functionals used in the definition of projection depth; see Zuo \cite{Zuo:2003p258}.
\end{appendix}


\section*{Acknowledgements}
Davy Paindaveine's research is supported by an A.R.C. contract
from the Communaut\'e Fran\c{c}aise de Belgique and by the IAP
research network grant nr. P7/06 of the Belgian government
(Belgian Science Policy). Germain Van Bever's research is supported
through a Mandat d'Aspirant FNRS (Fonds National pour la Recherche
Scientifique), Communaut\'e Fran\c{c}aise de Belgique. The authors are
grateful to an anonymous referee and the Editor for their careful
reading and insightful comments that allowed to improve the original manuscript.


%

\printhistory

\end{document}